\documentclass[12pt]{article}
\usepackage{amsmath}
\usepackage{amssymb}
\usepackage{mathrsfs}
\usepackage{times}
\usepackage{epsfig,graphics}
\usepackage{color}

\usepackage{amsthm}

\begin{document}
\pagestyle{plain}
\hsize = 6.5 in
\vsize = 8.5 in
\hoffset = -0.5 in
\voffset = -0.5 in
\baselineskip = 0.29 in
\global\long\def\rd{{\rm d}}
\def\vb {{\bf b}}
\def\vf {{\bf f}}
\def\vj {{\bf j}}
\def\vn {{\bf n}}
\def\vP{{\bf p}}
\def\vp{{\bf p}}
\def\vq{{\bf q}}
\def\vx{{\bf x}}
\def\vy{{\bf y}}
\def\vu{{\bf u}}
\def\vv{{\bf v}}
\def\vw{{\bf w}}
\def\vz {{\bf z}}
\def\A{{\bf A}}
\def\B{{\bf B}}
\def\mC{{\bf C}}
\def\mD{{\bf D}}
\def\F{{\bf F}}
\def\mG{{\bf G}}
\def\I{{\bf I}}
\def\vJ{{\bf J}}
\def\mM{{\bf M}}
\def\mN{{\bf N}}
\def\mS{{\bf S}}
\def\mQ{{\bf Q}}
\def\mR{{\bf R}}
\def\mT{{\bf T}}
\def\mU{{\bf U}}
\def\mV{{\bf V}}
\def\vX{{\bf X}}
\def\wtA{\widetilde{A}}
\def\mLa{\mbox{\boldmath$\Lambda$}}
\def\mGa{\mbox{\boldmath$\Gamma$}}
\def\mPhi{\mbox{\boldmath$\Phi$}}
\def\mPi{\mbox{\boldmath$\Pi$}}
\def\mXi{\mbox{\boldmath$\Xi$}}
\def\vg {\mbox{\boldmath$\gamma$}}
\def\vpi{\mbox{\boldmath$\pi$}}
\def\vphi{\mbox{\boldmath$\phi$}}
\def\vnu{\mbox{\boldmath$\nu$}}
\def\vkappa{\mbox{\boldmath$\kappa$}}
\def\vmu{\mbox{\boldmath$\mu$}}
\def\vxi{\mbox{\boldmath$\xi$}}
\def\vzeta{\mbox{\boldmath$\zeta$}}
\def\epr{\mbox{epr}}
\def\( {\left( }
\def\) {\right) }
\def\dbar{\ {\mathchar'26\mkern-12mu \rd}}

\title{The Statistical Foundation of Entropy in Extended Irreversible
Thermodynamics}

\author{Liu Hong\\[5pt]
Zhou Pei-Yuan Center for Applied Mathematics,\\
Tsinghua University, Beijing, 100084, P.R.C.\\
Email: zcamhl@tsinghua.edu.cn\\[10pt]
and\\[10pt]
Hong Qian\\[5pt]
Department of Applied Mathematics,\\
University of Washington, Seattle, WA 98195-3925, U.S.A.\\
Email: hqian@uw.edu
}


\maketitle

\tableofcontents

\begin{abstract}
In the theory of \textit{extended irreversible thermodynamics} (EIT), the flux-dependent entropy function plays a key role and has a fundamental distinction from the usual flux-independent entropy function adopted by \textit{classical irreversible thermodynamics} (CIT). However, its existence, as a prerequisite for EIT, and its statistical origin have never been justified. In this work, by studying the macroscopic limit of an $\epsilon$-dependent Langevin dynamics, which admits a large deviations (LD) principle, we show that the stationary LD rate functions of probability density $p_{\epsilon}(x,t)$ and joint probability density $p_{\epsilon}(x,\dot{x},t)$ actually turn out to be the desired flux-independent entropy function in CIT and flux-dependent entropy function in EIT respectively. The difference of the two entropy functions is determined by the time resolution for Brownian motions times a Lagrangian, the latter arises from the LD Hamilton-Jacobi equation and can be used for constructing conserved Lagrangian/Hamiltonian dynamics. 
\end{abstract}

\section{Introduction}

Classical irreversible, nonequilibrium thermodynamics for
macroscopic systems championed by the so called
Belgian-Dutch school, developed by Onsager, Meixner, Prigogine,
and many other authors, is based on the {\em local equilibrium hypothesis} \cite{degroot}.  The supposition guarantees
the existence of an entropy function $S(\vu)$ of the
macroscopic state variable $\vu$, which itself can be
a function of space $\vx$ and time $t$ in a system with
irreversible transport \cite{chapman-cowling}.
To go beyond the local equilibrium hypothesis,
{\em extended irreversible thermodynamics} assumes the existence of a new type of entropy functions
$\mathcal{S}(\vu,\vq)$ where variable $\vq$ is a flux
that represents the rates of transport processes \cite{Muller, jou}.
In classical thermodynamics, the very existence of a
``thermodynamic potential function'', as a principle, is
sufficient for deriving a collection of mathematical
relations that encompass the physics of thermodynamics.

	To provide the abstractly introduced entropy function
with a mechanistic basis, Helmholtz and Boltzmann
advanced the mechanical theory of heat
which firmly established that the concept of entropy in
thermodynamics has a statistical foundation in terms of
the dynamics of the constituents of a macroscopic system.
They were able to mathematically
derive the Gibbs' equation $\rd E = T\rd S -p\rd V$ for
mechanical systems in thermodynamic equilibrium based
on  (i) identifying a thermodynamic state as an entire
level set of a Hamiltonian function $H(\vx,\vp)$; and
(ii) Boltzmann's entropy $S(E)=k_B\ln\Omega(E)$,
where $\Omega(E)$ is the Lebesgue volume of
$\{(\vx,\vp)| H(\vx,\vp)\le E\}$ \cite{gallavotti-book}.

	In recent years, replacing the deterministic
Hamiltonian description by a stochastic Markov dynamics
and identifying the Gibbs-Shannon entropy as a mesoscopic
counterpart of entropy in a system with fluctuations, a rather
complete nonequilibrium thermodynamics in a state
space has been formulated \cite{qkkb}.
This theory exhibits four novel features:
(i) It represents all transport phenomena universally as the
probabilistic flux in the state space; then entropy production $=$
entropic force $\times$ probabilistic flux.
(ii) It removes the need for the local equilibrium assumption;
in fact it shows that the assumption is only a part of developing
Markovian models for real world processes as engineering.
(iii) It proves an ``law of entropy
balance'' \cite{degroot} as a theorem, providing the notions of
entropy production and entropy exchange with a
stochastic dynamic representation.  (iv) If the Markov process
has detailed balance, then the entropy exchange becomes
the rate of a mean potential energy change.

	In the light of this development,
``what is the statistical foundation of the $\mathcal{S}(\vu,\vq)$
in EIT?'' In the present work, we extend the stochastic,
Markov formulation of irreversible thermodynamics to address this important question.  There should be no doubt that the statistical
foundation of the $\mathcal{S}(\vu,\vq)$ has to reside in a
stochastic dynamics of the constituents of a mesoscopic system.

\section{Mesoscopic stochastic dynamics and its macroscopic limit}

\subsection{Mesoscopic stochastic dynamics}

	By {\em mesoscopic}, we mean a dynamic description of
a system in terms of a stochastic mathematical represenation,
with either discrete or continuous state space and time.
In the present work, we consider only the continuous time.
We give the general formalism in a continuous state space
$\mathbb{R}^n$, which in fact covers discrete, integer-valued
$\mathbb{Z}^n$ using Dirac-$\delta$ function.  To clearly illustrate
our ideas, more
involved mathematical derivations in the second part of the
paper, however, are carried out in terms of a discrete state space.

To be specific, let us consider a continuous-time, stochastic Markov dynamics in
a state space $\mathfrak{S}$, which are completely specified
by two mathematical objects: A probability
distribution $p(\vx,0)$, as an initial condition, and a transition rate
function for the probability $T(\vx,t+\Delta t|\vx',t)$.
All information concerning transport processes in the
state space $\mathfrak{S}$ is coded in the function
$T: \mathfrak{S}\to\mathfrak{S}$, and
\begin{equation}
    p(\vx,t+\Delta t) = \int_{\mathfrak{S}} T(\vx,t+\Delta t|\vx',t) p(\vx',t)
           \rd\vx'.
\label{cke}
\end{equation}
Eq. \ref{cke} is known as Chapman-Kolmogorov equation.
It is the foundational equation for Markov dynamics.
In terms of this mathematical representation, Gibbs entropy
in statistical thermodynamics has been identified
as a functional of the $p(\vx,t)$:\footnote{For finite state space, this functional is also known as
Shannon entropy, which emerges from the asymptotic behavior
of the frequency distribution of $n$ identical, independently
distributed (i.i.d.) uniform random variables, as $n\to\infty$.}
\begin{equation}
              S^{\text{CIT}}[p(\vx,t)] = -\int_{\mathfrak{S}}
            p(\vx,t)\ln p(\vx,t) \rd\vx.
\label{gibbs-e}
\end{equation}
A rather complete CIT, without the local equilibrium hypothesis, has been developed based on Eq. (\ref{gibbs-e}) \cite{qkkb,gq-pre-10}.
One significant success of this theory is the unification of
discrete stochastic chemical kinetics and Gibbsian equilibrium
chemical thermodynamics and the extension of the latter to
open, living biochemical systems \cite{ge-qian-16,ge-qian-17}.

	How does the $\vq$ variable enter this stochastic
formalism?  Certainly all information concerning $\vq$
is contained in the $T$.  But it cannot be the rate
of transition probability {\em per se} since $\vq$ is
necessarily zero in an equilibrium.  One naturally
considers the ``net probability flux'' from
$\vx\to\vx'$
\begin{equation}
		  J(\vx',t+\Delta t|\vx,t) = p(\vx,t)T(\vx',t+\Delta t|\vx,t)-
              p(\vx',t)T(\vx,t+\Delta t|\vx',t),
\label{net-flux}
\end{equation}
which is zero if and only if a stochastic dynamical system
reaches equilibrium state with detailed balance.  For a
mescopic system, thus conceptually one expects the EIT
entropy is a function of both $p(\vx,t)$ and
$J(\vx',t+\Delta t|\vx,t)$, the two key characteristics of
a nonequilibrium system \cite{wangjin}.  The
$p(\vx,t)T(\vx',t+\Delta t|\vx,t)$ is called the one-way
flux from $\vx$ to $\vx'$, and the $J$ in (\ref{net-flux}) is
called the net flux from states $\vx$ to $\vx'$ \cite{hill-book}.

To address this issue, let us consider a stochastic process $\vx(t)$ given by the Langevin dynamics
\begin{align}
\rd\vx(t)=\vb(\vx)\rd t+\sqrt{2\epsilon \mD(\vx)}
 \rd\B(t),
\end{align}
with drift $\vb(\vx)$ and the diffusion coefficient $\mD(\vx)$ is symmetric and positive definite. $\epsilon\ll1$ is a small parameter indicating the level of stochasticity. As $\epsilon\rightarrow0$, the Langevin dynamics approaches to a deterministic dynamics $d\vx/dt=\vb(\vx)$. According to It\^{o}'s calculus, it is well-known that the instantaneous
probability density function $p_{\epsilon}(\vx,t)$ and
transition probability $T_{\epsilon}(\vx,t|\vx',t')$ both follow Kolmogorov forward equations
\begin{eqnarray}
\frac{\partial p_{\epsilon}(\vx,t)}{\partial t} &=&
\frac{\partial}{\partial \vx}\cdot\bigg[\epsilon \frac{\partial}{\partial \vx}\mD(\vx) p_{\epsilon}(\vx,t)-\vb(\vx)p_{\epsilon}(\vx,t)\bigg],
\label{e4p}\\
\frac{\partial T_{\epsilon}(\vx,t|\vx',t')}{\partial t}
 &=&
\frac{\partial}{\partial \vx}\cdot\bigg[\epsilon \frac{\partial}{\partial \vx}\mD(\vx)T_{\epsilon}(\vx,t|\vx',t')-\vb(x)T_{\epsilon}(\vx,t|\vx',t')\bigg].
\label{e4T}
\end{eqnarray}

\subsection{LDRF and classical irreversible thermodynamics}
The large deviations theory
supports a WKB ansatz, $p_{\epsilon}(\vx,t)=\exp[-\varphi(\vx,t)/\epsilon]$, based on which one
finds the large deviations rate function (LDRF) $\varphi(\vx,t)$ satisfies
a Hamilton-Jacobi equation (HJE)
\begin{equation}
  \frac{\partial \varphi(\vx,t)}{\partial t}=-\bigg[\frac{\partial\varphi(\vx,t)}{\partial \vx}\bigg]^T\mD(\vx)\frac{\partial\varphi(\vx,t)}{\partial \vx}-\bigg[\frac{\partial\varphi(\vx,t)}{\partial \vx}\bigg]^T\vb(\vx).
\end{equation}
It allows the introduction of a Hamiltonian function
\begin{align}
H(\vx,\vy)=\vy^T\mD(\vx)\vy+\vy^T\vb(\vx),
\end{align}
where $\vy=\partial\varphi(\vx,t)/\partial \vx$, and the corresponding Hamiltonian dynamics
\begin{align}
&\frac{d\vx}{dt}=\frac{\partial H(\vx,\vy)}{\partial \vy}=2\mD(\vx)\vy+\vb(\vx),\\
&\frac{d\vy}{dt}=-\frac{\partial H(\vx,\vy)}{\partial \vx}=-\vy^T\frac{\partial\mD(\vx)}{\partial \vx}\vy-\vy^T\frac{\partial\vb(\vx)}{\partial \vx}.
\end{align}
The Hamiltonian dynamics is a generalization of the deterministic dynamics $d\vx/dt=\vb(\vx)$, in which $\vy$ can be regarded as fluctuations in ``a momentum space''.  If
$\vy(0)=0$, then $\vy(t)=0$ and $\vx(t)$ follows the
$d\vx/dt=\vb(\vx)$. A very dramatic feature of this
generalization is the ``conservative nature''
of $(\vx,\vy)(t)$ dynamics.

If a diffusion process satisfies $\vb(\vx)
=-\mD(\vx)\nabla\varphi^{eq}(\vx)$, then it is sufficient and necessary that the diffusion
is non-driven.  The emergent Hamiltonian for a non-driven
stochastic system can be transformed, via a {\em canonical transformation}, into a form which is an even function of the
momentum variable $\vp$, signifying time reversibility:
\begin{subequations}
\begin{eqnarray}
	H(\vx,\vy) &=& \vy^T\mD(\vx)\vy + \vy^T\vb(\vx),
\nonumber\\
	&=& \vy^T\mD(\vx)\big(\vy-\nabla\varphi^{eq}(\vx)\big)
 = \vp^T\mD(\vq)\vp+ V(\vq) = \tilde{H}(\vq,\vp),
\label{newton-h}
\end{eqnarray}
in which $\vq=\vx$ and
\begin{equation}
       \vp = \vy-\frac{1}{2}\nabla\varphi^{eq}(\vx), \
       V(\vq) = -\frac{1}{4}[\nabla\varphi^{eq}(\vq)]^T\mD(\vq)\nabla\varphi^{eq}(\vq).
\end{equation}
\end{subequations}
Indeed, the Hamiltonian in (\ref{newton-h}) has the Newtonian
expression with a separation of a kinetic energy and a potential
energy.  The matrix $\mD(\vq)$ in the kinetic energy represents
a curved space.

To show the transformation is canonical, we note
\begin{equation}
     \Big\{\vx,\vy,H(\vx,\vy)\Big\} \to
     \Big\{\vx,\tilde{\vy}=\vy+\vf(\vx), \tilde{H}(\vx,\tilde{\vy})=H\big(\vx,\tilde{\vy}-\vf(\vx)\big)
      \Big\},
\end{equation}
has
\begin{eqnarray}
	\frac{\rd\vx}{\rd t} &=& \left(\frac{\partial H}{\partial \vy}\right)_{\vx} =  \left(\frac{\partial\tilde{H}(\vx,\tilde{\vy})}{\partial \tilde{\vy}}\right)_{\vx},
\\
	 \frac{\rd \tilde{\vy}}{\rd t} &=&
            \frac{\rd\vy}{\rd t} +\vf'_{\vx}(\vx)
           \left(\frac{\rd\vx}{\rd t}\right)
\nonumber\\
	&=& -\left(\frac{\partial H}{\partial \vx}\right)_{\vy}+
           \left(\frac{\partial H}{\partial\vy}\right)_{\vx}\vf'_{\vx}(\vx)
   = -\left(\frac{\partial \tilde{H}(\vx,\tilde{\vy})}{\partial\vx}\right)_{\tilde{\vy}}.
\end{eqnarray}

On the other hand, if a diffusion process has $\mD^{-1}(\vx)\vb(\vx)$ not
being a gradient vector field, then it is easy to show that its corresponding Lagrangian equation
\begin{subequations}
\begin{eqnarray}
D^{-1}_{ij}(\vx)\ddot{x}_j&=& \frac{1}{2}\underbrace{ \frac{\partial\big[ b_j(\vx)D^{-1}_{jk}(\vx)b_k(\vx)-\dot{x}_j D^{-1}_{jk}(\vx)\dot{x}_k\big]}{\partial x_i}}_{\text{ potential force}}
\\[7pt]
&+&  \underbrace{ \dot{x}_k\frac{\partial D^{-1}_{ij}(\vx)b_j(\vx)}{\partial x_k}-\dot{x}_j\frac{\partial D^{-1}_{jk}(\vx)b_k(\vx)}{\partial x_i}}_{\text{ Lorentz force 1}}
\underbrace{-\dot{x}_k\frac{\partial D^{-1}_{ij}(\vx)}{\partial x_k}\dot{x}_j+\dot{x}_j\frac{\partial D^{-1}_{jk}(\vx)}{\partial x_i}\dot{x}_k. }_{\text{ Lorentz force 2}}  \hspace{1cm}
\end{eqnarray}
\end{subequations}
has two Lorentz magnetic force like terms \cite{gq-ijmpb}, since they makes no contributions to the work ($\dot{\vx}\times$ Lorentz force $=0$).
One may also look for alternative time irreversible extensions, which is a central topic in nonequilibrium thermodynamics, by examining the stationary large deviations rate function $\varphi^{ss}(\vx(t))$,
\begin{equation}
	\bigg[\frac{\partial\varphi^{ss}(\vx)}{\partial \vx}\bigg]^T \left\{ \mD(\vx)\frac{\partial\varphi^{ss}(\vx)}{\partial \vx}+\vb(\vx) \right\} = 0.
\label{shje}
\end{equation}
Eq. \ref{shje} reveals a decomposition of the vector
field $\vb(\vx)$:
\begin{equation}
    \vb(\vx) = -\mD(\vx)\nabla\varphi^{ss}(\vx)
       + \vg(\vx),
\label{dcompb}
\end{equation}
in which $\vg^T(\vx) \cdot \nabla\varphi^{ss}(\vx)=0$
for all $\vx$.   The stationary large deviations rate
function $\varphi^{ss}(\vx)$ can, and should be identified as
the free energy function in irreversible thermodynamics, as
illustrated below.

Classical irreversible thermodynamics as presented
by Onsager and others first suggested that any
non-driven systems spontaneously approches to a
equilibrium steady state.   This corresponds to
$\vg(\vx)=0$ in the stochastic dynamics.  Then one has
\begin{equation}
      \frac{\rd\vx(t)}{\rd t} = -\mD(\vx)\nabla\varphi^{ss}(\vx).
\label{graddy}
\end{equation}
This is precisely what has been expected from and discussed
in CIT.  In fact,
\begin{equation}
   \frac{\rd}{\rd t}\varphi^{ss}\big(\vx(t)\big)
   = \left[\frac{\partial\varphi^{ss}(\vx)}{\partial\vx}\right]^T \frac{\rd\vx(t)}{\rd t}
   = -\vb^T(\vx)\mD^{-1}(\vx)\vb(\vx) \le 0,
\end{equation}
which means $\varphi^{ss}$ is the relative entropy for CIT, since it is also positive and convex as a fundamental mathematical property of LDRF.
The theory of CIT particularly recognizes a geometric
interpretation of $\mD^{-1}(\vx)$: It provides an
appropriate metric in the tangent space of $\vx$, to which
$\vb(\vx)$ belongs.

	More generally without detailed balance, based on
(\ref{dcompb}) one still has
\begin{equation}
    \frac{\rd}{\rd t}\varphi^{ss}\big(\vx(t)\big)
 = \left[\frac{\partial\varphi^{ss}(\vx)}{\partial\vx}\right]^T
 \vb(\vx) = -\left[\frac{\partial\varphi^{ss}(\vx)}{\partial\vx}\right]^T
 \mD(\vx)\left[\frac{\partial\varphi^{ss}(\vx)}{\partial\vx}\right]
     \le 0.
 \label{flux-inde}
\end{equation}
In fact,
\begin{equation}
  \frac{\rd}{\rd t}\varphi^{ss}\big(\vx(t)\big)
    =  -\vb^T(\vx)\mD^{-1}(\vx)\vb(\vx) +
       \vg^T(\vx)\mD^{-1}(\vx)\vg(\vx),
\label{febe}
\end{equation}
which implies a Pythagorean relation among three
entropy productions \cite{qian-arXiv}:
\begin{equation}
	\underbrace{ \vb^T(\vx)\mD^{-1}(\vx)\vb(\vx)
}_{\text{total entropy production}} =  \underbrace{ \big[\mD(\vx)\nabla\varphi^{ss}\big]^T
       \mD^{-1}(\vx)\big[\mD(\vx)\nabla\varphi^{ss}\big] }_{\text{free energy dissipation}}
           + \underbrace{ \vg^T(\vx)\mD^{-1}(\vx)\vg(\vx) }_{\text{house-keeping heat}}.
\label{epdcp}
\end{equation}
The two terms on the rhs of (\ref{epdcp}) have been
identified as {\em Boltzmann's thesis} and {\em Prigogine's
thesis} of irreversibility \cite{qian-book-chapter}.
Boltzmann's thesis focuses on
transient relaxation dynamics that approaches to an equilibrium
in a non-driven system, and Prigogine's idea that articulates
driven phenomena that can exist even in a stationary
state \cite{gq-pre-10}. In stochastic thermodynamics, these two origins are represented
by {\em free enegy dissipation} and {\em house-keeping
heat}, as two distinct parts of the total entropy production. The
house-keeping heat has a dual interpretation: as an
external driving force to an overdamped thermodynamics
or as the inertia effect in a conservative
dynamics \cite{qian-pla-2014}.  The latter interpretation,
as we show below, can be further developed
in terms of an internal conjugate variable.

	A remark is in order: Eq. \ref{febe} is a more legitimate
thermodynamic law than the entropy balance
equation \cite{degroot}:
\begin{equation}
                       \frac{\rd S}{\rd t}= \frac{\dbar_iS}{\rd t}
                +\frac{\dbar_e S}{\rd t},
\label{ebl}
\end{equation}
in which among the three terms, entropy change ($\rd S$),
entropy production ($\dbar_iS$), and entropy flux ($\dbar_eS$),
only the $(\dbar_i S/\rd t)$ has
a definite sign.  The $\varphi^{ss}$ on the
lhs of (\ref{febe}) is a free energy, and each one of the
three terms in (\ref{febe}) has a definit sign.  As it has been
known from equilibrium thermodynamics, free energy is the appropriate thermodynamic potential function of a non-isolated system; not entropy.

\section{Flux-dependent entropy and irreversible thermodynamics}

\subsection{Flux-dependent entropy function}

	For stochastic dynamics without detailed balance, Eq.
\ref{graddy} is no longer true.  In order to take the nonzero vector field $\vg(\vx)$ into consideration, one needs to study not only the state of a system, but also the fluxes between any two given states. Flux of a mesoscopic stochastic dynamics in $\mathfrak{S}$, as
given in (\ref{net-flux}), is completely determined by $p(\vx,t)$
and transition probability $T(\vx,t+\Delta t|\vx',t)$ defined in
(\ref{e4p}) and (\ref{e4T}).  For an infinitesimal $\Delta t$,
the transition probability for the diffusion process
in (\ref{e4T}) has the form
\begin{align}
&T_{\epsilon}(\vx,t+\Delta t|\vx',t)\nonumber\\
&=\frac{1}{
    \sqrt{(4\pi\epsilon(\Delta t))^n |\det\mD(\vx')|} }
          \exp\left[-\frac{[\vx-\vx'-\textbf{\~{b}}(\vx')\Delta t]^T
               \mD^{-1}(\vx')[\vx-\vx'-\textbf{\~{b}}(\vx')\Delta t]}{4\epsilon(\Delta t)}\right].
\label{T-gauss}
\end{align}
Then Hill's net probability flux in (\ref{net-flux}) $(\Delta\vx)^{-1}J(\vx+\Delta\vx,t+\Delta t|\vx,t)$ becomes
\begin{eqnarray}
  && \lim_{\Delta\vx\to 0}
   \frac{1}{|\Delta\vx|}\big[ p_{\epsilon}(\vx,t) T_{\epsilon}
         (\vx+\Delta\vx,t+\Delta t |\vx,t)
            - p_{\epsilon}(\vx+\Delta\vx,t)T_{\epsilon}
         (\vx,t+\Delta t |\vx+\Delta\vx,t) \big]
\nonumber\\
	&=&
[\epsilon\mD(\vx)]^{-1}\vb(\vx)p_{\epsilon}(\vx,t)-\nabla p_{\epsilon}(\vx,t).
\end{eqnarray}
The probability flux in diffusion theory,
$\vJ_{\epsilon}(\vx,t)\equiv \vb(\vx)p_{\epsilon}(\vx,t)-\epsilon\mD(\vx)\nabla p_{\epsilon}(\vx,t)$ is actually $\epsilon\mD(\vx)$ $\times$ Hill's net flux. This result also reveals that while the mesoscopic flux $J(\vx+\Delta\vx,t+\Delta t|\vx,t)$ is completely determined once $p_{\epsilon}(\vx,t)$ and the transition
probability $T_{\epsilon}(\vx,t|\vx',t')$ are know, in the macroscopic limit, the transport
flux $J_{\epsilon}(\vx,t)$ is not determined by $\vx(t)$ and vector field
$\vb(\vx)$.

In the macroscopic limit as $\epsilon\to 0$, it can be shown that
\begin{subequations}
\begin{eqnarray}
      p_{\epsilon}(\vx,t) &\to& \delta\big(\vx-\tilde{\vz}(t)\big),
\\
	- \epsilon\ln p_{\epsilon}(\vx,t) &\to& \varphi(\vx,t),
\\
	-\epsilon\ln T_{\epsilon}(\vx+\Delta\vx,t+\Delta t|\vx, t) &\to&
            L(\vx,\dot{\vx})\Delta t,
\\
	- \epsilon\ln p_{\epsilon}^{ss}(\vx) &\to& \varphi^{ss}(\vx),
\\
    \frac{\vJ^{ss}_{\epsilon}(\vx)}{p^{ss}_{\epsilon}(\vx) } &\to&
       \vg(\vx),
\end{eqnarray}
\end{subequations}
in which $\tilde{\vz}(t)$ is the solution to $\dot{\vz}=\vb(\vz)$, $\vg(\vx)$ is defined in (\ref{dcompb}). The Lagrangian
$L(\vx, \dot{\vx})=\frac{1}{4}[\dot{\vx}-\vb(\vx)]^T
\mD^{-1}(\vx)[\dot{\vx}-\vb(\vx)]$, with $\dot{\vx}=\Delta\vx/\Delta t$ as the large deviations rate function for the
transition probability over infinitesimal $\Delta t$. And,
\begin{equation}
  \frac{ \vJ_{\epsilon}(\vx,t) }{p_{\epsilon}(\vx,t)}
   = \vb(\vx) -\epsilon\mD(\vx)\nabla \ln  p_{\epsilon}(\vx,t)
   \to \vb(\vx)+\mD(\vx)\nabla\varphi(\vx,t).
\end{equation}

With respect to the probability density $p_{\epsilon}(\vx,t)$
and transition probability $T_{\epsilon}(\vx,t+\Delta t|\vx',t)$, one natural choice of the flux-dependent entropy function is
\begin{eqnarray}
S^{meso-EIT}(t;\Delta t)=-\int_{\mathfrak{S}}p_{\epsilon}(\vx,t)T_{\epsilon}(\vx',t+\Delta t|\vx,t)\ln[p_{\epsilon}(\vx,t)T_{\epsilon}(\vx',t+\Delta t|\vx,t)]d\vx d\vx'.
\end{eqnarray}
Then we have
\begin{eqnarray}
 && S^{meso-EIT}(t;\Delta t)
\nonumber\\
&=& S^{meso-CIT}(t;\Delta t)-\int_{\mathfrak{S}}p_{\epsilon}(\vx,t)T_{\epsilon}(\vx',t+\Delta t|\vx,t)\ln T_{\epsilon}(\vx',t+\Delta t|\vx,t)d\vx d\vx'.
\end{eqnarray}
Considering the Gaussian-form solution of $T_{\epsilon}(\vx',t+\Delta t|\vx,t)$ given in (\ref{T-gauss}), it is easy to show that the
difference between $S^{\text{meso-EIT}}$
and $S^{\text{meso-CIT}}$ is a function of $\Delta t$, which
is expected to tending zero as $\Delta t\to 0$.\footnote{Note that
for a continuous distribution, the mathematics of
\[
    \lim_{\epsilon\to 0} \int_{\mathbb{R}} p_{\epsilon}(x)\ln p_{\epsilon}(x)
\rd x
\]
where $p_{\epsilon}(x)=\delta(x)$ as $\epsilon\to 0$, is not
necessarily zero!  An example is the Gaussian distribution
with variance $\epsilon$:
\begin{equation*}
	 -\int_{\mathbb{R}} p_{\epsilon}(x)\ln p_{\epsilon}(x) \rd x	
  = -\int_{\mathbb{R}} p_{\epsilon}(x)
      \left[ -\frac{x^2}{2\epsilon} -\frac{1}{2}\ln (2\pi\epsilon)  \right] \rd x	 = \frac{1}{2}+\frac{1}{2}\ln (2\pi\epsilon).
\end{equation*}
It is not zero; it does not even converge as $\epsilon\to 0$.
This is in sharp contrast to a discreate distribution, which
has  $\sum_{i} \delta_{i,i_0}\ln\delta_{i,i_0} = 0$.}
Following the definition, we further have
\begin{align}
&\frac{dS^{meso-CIT}}{dt}=-\int_{\mathfrak{S}}\frac{\partial p_{\epsilon}(\vx,t)}{\partial t}\ln p_{\epsilon}(\vx,t)d\vx=-\int_{\mathfrak{S}}\frac{\partial}{\partial \vx}\bigg[\epsilon \mD(\vx)\frac{\partial p_{\epsilon}}{\partial\vx}-\vb(\vx)p_{\epsilon}\bigg]\ln p_{\epsilon}dx\nonumber\\
=&\int_{\mathfrak{S}}\bigg[\epsilon \mD(\vx)\frac{\partial p_{\epsilon}}{\partial\vx}-\vb(\vx)p_{\epsilon}\bigg]\frac{\partial}{\partial \vx}(\ln p_{\epsilon})dx=-\int_{\mathfrak{S}}\vb(\vx)\frac{\partial p_{\epsilon}}{\partial \vx}dx+\int_{\mathfrak{S}}\frac{\partial p_{\epsilon}}{\partial \vx}[\epsilon \mD(\vx)]\frac{\partial p_{\epsilon}}{\partial\vx}dx.
\end{align}
The last two terms represent entropy flux and entropy production rate respectively. Meanwhile,
\begin{align}
&\frac{dS^{meso-EIT}}{dt}-\frac{dS^{meso-CIT}}{dt}\nonumber\\
=&-\int_{\mathfrak{S}}\bigg[\frac{\partial}{\partial t}p_{\epsilon}(\vx,t)T_{\epsilon}(\vx',t+\Delta t|\vx,t)+p_{\epsilon}(\vx,t)\frac{\partial}{\partial t}T_{\epsilon}(\vx',t+\Delta t|\vx,t)\bigg]\ln T_{\epsilon}(\vx',t+\Delta t|\vx,t)d\vx d\vx'.\nonumber
\end{align}

%
%
%
%

\subsection{Lagrangian function and conditional probabilities}

	In the theory of large deviations, $\epsilon$ stands for the
level of stochasticity.  The large deviations principle then states $e^{-\varphi(\vx,t)/\epsilon}$
as the probability density function of $\vx_{\epsilon}(t)$, with
the rate function $\varphi(\vx,t)$ given as
\begin{equation}
   \varphi(\vz,t) = \min_{\tiny\begin{array}{c}
          \vx(s) \\ \vx(0)=\vx_0 \\ \vx(t)=\vz \end{array} }
     \int_0^{t}
                L\big[\vx(s),\dot{\vx}(s)\big]\rd s.
\label{equation021}
\end{equation}
Let us particularly consider $t$ to $t+2\Delta t$ with
a very small $\Delta t$.  Then,
\begin{eqnarray*}
	  && \min_{\tiny\begin{array}{c}
          \vx(s), t\le s\le t+2\Delta t \\ \vx(t)=\vx_0 \\ \vx(t+2\Delta t)=\vx_1 \end{array} }
        \int_{t}^{t+2\Delta t} L\big[\vx(s),\dot{\vx}(s)\big]\rd s
\\
      &=&  \min_{\tiny\begin{array}{c}
          \vx(t+\Delta t)  \end{array} }      \Delta t
   \Big\{ L\left[\vx_0,\tfrac{\vx(t+\Delta t)-\vx_0}{\Delta t}\right]  +
 L\left[\vx(t+\Delta t),\tfrac{\vx_1-\vx(t+\Delta t)}{\Delta t}\right]
             \Big\}
\\
      &=&    \Delta t\Big\{ L\left[\vx_0,\tfrac{\vx^*-\vx_0}{\Delta t}\right]  +
 L\left[\vx^*,\tfrac{\vx_1-\vx^*}{\Delta t}\right]
             \Big\},
\end{eqnarray*}
in which
\begin{equation*}
  \left(\frac{\partial L(\vx,\dot{\vx})}{\partial\vx}\right)_{\vx^*,\tfrac{\vx_1-\vx^*}{\Delta t}}+
  \frac{1}{\Delta t} \left[ \left(\frac{\partial L(\vx,\dot{\vx})}{\partial \dot{\vx}} \right)_{\vx_0,\tfrac{\vx^*-\vx_0}{\Delta t}} - \left(\frac{\partial L(\vx,\dot{\vx})}{\partial \dot{\vx}}
   \right)_{\vx^*,\tfrac{\vx_1-\vx^*}{\Delta t}}\right]
  = 0.
\end{equation*}
In the limit of $\Delta t\to 0$, this recovers
the Euler-Lagrange equation,
\begin{equation}
 \left(\frac{\partial L(\vx,\dot{\vx})}{\partial\vx}\right)
   -\frac{\rd}{\rd t}\left(\frac{\partial L(\vx,\dot{\vx})}{\partial \dot{\vx}}
   \right) = 0.
\end{equation}

For a diffusion process, its Kolmogorov forward equation gives us at
$\vx$ the transition probability for a $\Delta\vx$ during a
sufficiently small $\Delta t$, the conditional probability
density for the $\dot{\vx}$, thus,
is a Gaussian distribution with mean $\vb(\vx)$ and
covariance matrix $(2\epsilon/\Delta t)\mD(\vx)$ \cite{gq-ijmpb}:
\begin{eqnarray}
   p_{\epsilon}(\dot{\vx}|\vx;\Delta t)  &=&
  \frac{1}{\sqrt{ (4\pi\epsilon/\Delta t)^n\det\mD(\vx) }} \exp\left\{-\frac{\Delta t[\dot{\vx}-\vb(\vx)]^T\mD^{-1}(\vx)[\dot{\vx}-\vb(\vx)] }{4\epsilon}\right\}
\nonumber\\
	    &=&  \frac{1}{\sqrt{ (4\pi\epsilon/\Delta t)^n\det\mD(\vx) }} e^{-\frac{\Delta t}{\epsilon}L[\vx(t),\dot{\vx}(t)]},
\label{eqn023}
\end{eqnarray}
in which
\begin{equation}
        L\big[\vx(s),\dot{\vx}(s)\big]= \frac{1}{4}
        \big[\dot{\vx}-\vb(\vx) \big]^T \mD^{-1}(\vx)
             \big[\dot{\vx}-\vb(\vx) \big].
\label{equation024}
\end{equation}
The probabilistic significance of $e^{-\frac{\Delta t}{\epsilon}L(\vx,\dot{\vx})}$ is to provide the probability of the
$\dot{\vx}\equiv\frac{\Delta\vx}{\Delta t}$,
conditioned at $\vx$, with the ``time resolution'' $\Delta t$. Mathematically, this means we consider $\vx_{\epsilon}(t)$
in the context of ``certain smooth functions'' while strictly
speaking, according to It\^{o}, $\vx_{\epsilon}(t)$
is non-differentiable!

Now noting the relation between
$\vy$ and $\dot{\vx}=\vb(\vx)+2\mD(\vx)\vy,$
the conditional probability density for the conjugate variable
$\vy$, or momentum, is
\begin{equation}
 p_{\epsilon}(\vy|\vx;\Delta t)
  = \frac{1}{\sqrt{(\pi\epsilon/\Delta t)^n\det\mD^{-1}(\vx) }}
       \exp\left\{-\frac{\Delta t\vy^T\mD(\vx)\vy }{\epsilon}\right\},
\label{eqn025}
\end{equation}
which is again a Gaussian distribution, with zero mean and
covariance matrix $\epsilon/(2\Delta t)\mD^{-1}(\vx)$.
It is noted that the covariance matrices for $\dot{\vx}$ and $\vy$
are proportional to $\mD(\vx)$ and $\mD^{-1}(\vx)$ respectively.

It is easy to verify the familiar relationship between Lagrangian
$L(\vx,\dot{\vx})$ in (\ref{equation024}) and Hamiltonian
function $H(\vx,\vy)=\vy^T\mD(\vx)\vy+\vy^T\vb(\vx)$:
\begin{eqnarray*}
	 \vy &=& \left(\frac{\partial L(\vx,\dot{\vx})}{\partial \dot{\vx}}
       \right)_{\vx} =
     \frac{1}{2}\mD^{-1}(\vx)\big[\dot{\vx}-\vb(\vx) \big],
\\
	\dot{\vx} &=& \left( \frac{\partial H(\vx,\vy)}{\partial\vy}
          \right)_{\vx} = 2\mD(\vx)\vy + \vb(\vx),
\\
     L\big(\vx,\dot{\vx}\big) &=&  \dot{\vx}^T\vy
          - H(\vx,\vy).
\end{eqnarray*}

	From the conditional probability in (\ref{eqn023}) and
(\ref{eqn025}), we have the joint probability density function
for $\vx$ and $\dot{\vx}$:
\begin{equation}
	p_{\epsilon}(\vx,\dot{\vx},t;\Delta t) =
   A_1(\epsilon,t,\Delta t)\exp\left\{ -\frac{1}{\epsilon}
             \Big[ \varphi(\vx,t)+\Delta tL(\vx,\dot{\vx})\Big]\right\},
\end{equation}
in which $A_1(\epsilon,t,\Delta t)$ is a normalization factor.
Similarly, the joint probability for $\vx$ and $\vy$:
\begin{equation}
	p_{\epsilon}(\vx,\vy,t;\Delta t) =
    A_2(\epsilon,t,\Delta t)\exp\left\{ -\frac{1}{\epsilon}
             \Big[ \varphi(\vx,t)+\Delta t\vy^T\mD(\vx)\vy\Big]\right\}.
\label{eqn27}
\end{equation}
It should be noted that $\Delta t$ stands for the time resolution required for the existence of a normal ``smooth'' diffusion process. Roughly speaking, as the Brownian motion is non-differentiable, in order to properly define $\vx$ and $\dot{\vx}$ in the context of ``certain smooth functions'', we need to coarse grain the time scale by looking at their averages over a ``microscopically sufficiently long yet macroscopically sufficiently short'' (due to CIT) time resolution $\Delta t$. The shorter $\Delta t$ is, the larger $\dot{\vx}$ (or $\vy$) will be, as a manifestation of certain uncertainty principle we will address in detail later. In this sense, even though $\Delta t\ll1$, $\Delta t\vy^T\mD(\vx)\vy$ may still be comparable with $\varphi(\vx,t)$ and makes a non-negligible contribution to the joint probability.

\subsection{LDRF and extended irreversible thermodynamics}

To go beyond the so-called
{\em local equilibrium hypothesis}, the extended irreversible thermodynamics has been proposed by M\"{u}ller and Ruggeri \cite{Muller}, Jou, Casas-V\'{a}zquez and Lebon \cite{jou}, as a modification of classical irreversible thermodynamics. A major difference of the two theories lays on the choice of state variables. In CIT, only variables used in equilibrium thermodynamics are allowed, while in EIT nonequilibrium variables characterizing the fluxes of transport processes are adopted too. For example, in a EIT formulation of classical hydrodynamics, the fluid density $\rho$, velocity $v$, total energy $E$, stress tensor $P$ and heat flux $q$ are all taken as independent variables. While, in CIT the stress tensor $P$ and heat flux $q$ have to be treated as dependent variables, \textit{i.e.} $P=P(\rho, v, E)$ and $q=q(\rho, v, E)$. This difference is raised by the fact that only the first three variables appear in the description of equilibrium thermodynamics of fluids, while the latter two are not. Actually, $P$ and $q$ are fluxes relating to the transport of momentum and energy in a nonequilibrium process.

As EIT adopts an enlarged space of state variables, it shows a stronger power in dealing with nonequilibrium processes than CIT. A first non-trivial successful application of EIT is the derivation of Cattaneo's law for heat conduction, which solves the problem of infinite-speed propagation of thermal signals obtained from the Fourier's law. Later, EIT has been applied to a rich phenomenology in heat transport, second sound in solids, ultrasound propagation or generalized hydrodynamics, \textit{ etc.} \cite{jou}. Despite its great success, the origin of flux-dependent entropy function in EIT has never been clarified. Interestingly, as we have shown above, the large derivation function obtained from the limit process of a mesoscopic stochastic dynamics turns out to be the entropy function for CIT-like modeling theories. Therefore, we would like to see a possible emergence of a set of EIT-like theories in our stochastic framework too.

To make this point clear, we look for large derivation functions as a function of both state variable $\vx$ and its flux in accordance with EIT. Obviously, the conditional probability in (\ref{eqn27}) meets our requirement, \textit{i.e.}
\begin{align}
\varphi(\vx,\vy,t;\Delta t)=-\lim_{\epsilon\rightarrow0}\epsilon\ln[p_{\epsilon}(\vx,\vy,t;\Delta t)]=\varphi(\vx,t)+\vy^T[\Delta t\mD(\vx)]\vy,
\end{align}
in which $\vy=\frac{1}{2}\mD^{-1}(\vx)[\dot{\vx}-\vb(\vx)]$. $\varphi(\vx,\vy,t;\Delta t)$ can be regarded as a Level 1.5 LDRF, since the ordinary Level 1 LDRF
\begin{align}
\varphi(\vx,t)=\min_{\vy}\varphi(\vx,\vy,t;\Delta t)
\end{align}
can be obtained by the contraction principle. Interestingly, it is easy to see that the minimum in the above formula is reached at $\vy=0$ or $d\vx/dt=\vb(\vx)$, the determinist dynamics when $\epsilon=0$.

The stationary large derivation rate function
\begin{align}
\varphi^{ss}(\vx,\vy;\Delta t)=\lim_{t\rightarrow\infty} \varphi(\vx,\dot{\vx},t;\Delta t)=\varphi^{ss}(\vx)+\vy^T[\Delta t\mD(\vx)]\vy
\end{align}
actually provides the statistical foundation of the flux-dependent entropy function used in the extended irreversible thermodynamics. Its full time derivative obeys the entropy balance law,
\begin{align}
&\frac{d\varphi^{ss}(\vx,\vy,\Delta t)}{dt}=\bigg[\frac{\partial\varphi^{ss}(\vx,\vy,\Delta t)}{\partial \vx}\bigg]^T\frac{d\vx}{dt}+\bigg[\frac{\partial\varphi^{ss}(\vx,\vy,\Delta t)}{\partial \vy}\bigg]^T\frac{d\vy}{dt}\nonumber\\
=&\bigg[\frac{\partial\varphi^{ss}(\vx)}{\partial \vx}\bigg]^T\frac{d\vx}{dt}+\Delta t\bigg[\frac{\partial\big[\vy^T\mD(\vx)\vy\big]}{\partial \vx}\bigg]^T\frac{d\vx}{dt}+2 \vy^T[\Delta t\mD(\vx)]\frac{d\vy}{dt}\nonumber\\
=&\bigg\{\bigg[\frac{\partial\varphi^{ss}(\vx)}{\partial \vx}\bigg]^T+\Delta t\bigg[\frac{\partial\big[\vy^T\mD(\vx)\vy\big]}{\partial \vx}\bigg]^T\bigg\}\bigg[2\mD(\vx)\vy+\vb(\vx)\bigg]+2 \vy^T[\Delta t\mD(\vx)]\frac{d\vy}{dt}\nonumber\\
=&\frac{\partial}{\partial\vx}\cdot\bigg[\Delta t\big(\vy^T\mD(\vx)\vy\big)\big(2\mD(\vx)\vy+\vb(\vx)\big)\bigg]+\bigg[\frac{\partial\varphi^{ss}(\vx)}{\partial \vx}\bigg]^T\vb(\vx)\nonumber\\
&+2\Delta t\vy^T\mD(\vx)\bigg\{(\Delta t)^{-1}\frac{\partial\varphi^{ss}(\vx)}{\partial\vx}-\frac{1}{2}\frac{\partial}{\partial\vx}\cdot\bigg[2\mD(\vx)\vy+\vb(\vx)\bigg]\vy+\frac{d\vy}{dt}\bigg\},
\end{align}
by inserting the known relation $d\vx/dt=2\mD(\vx)\vy+\vb(\vx)$ and using integration by parts. In last step, the first term is recognized as the entropy flux. The next two terms are entropy production rates, which must be non-positive and equal to zero if and only if at the stationary state. Actually, it has already been shown in (\ref{flux-inde}) that $[\partial\varphi^{ss}(\vx)/\partial \vx]^T\vb(\vx)\leq0$ in accordance with classical irreversible thermodynamics, so that we only need to require
\begin{align}
\frac{d\vy}{dt}=-(\Delta t)^{-1}\frac{\partial\varphi^{ss}(\vx)}{\partial\vx}+\frac{1}{2}\frac{\partial}{\partial\vx}\cdot\bigg[2\mD(\vx)\vy+\vb(\vx)\bigg]\vy-\alpha(\vx,\vy)\vy
\end{align}
where $\alpha(\vx,\vy)\geq0$ is a non-negative function. In particular, at the stationary state when $d\vy/dt=0$, we arrive at the gradient dynamics,
\begin{align}
\vy\propto(\Delta t)^{-1/2}\frac{\partial\varphi^{ss}(\vx)}{\partial\vx},
\end{align}
which happens at the correct time scale of $(\Delta t)^{-1/2}$ for Brownian motions.

The global minimum of $\varphi^{ss}(\vx,\vy;\Delta t)$ is obtained at
\begin{equation}
    \vy^* = 0, \  \vx^* = \min_{\vx}\varphi^{ss}(\vx).
\end{equation}
For stochastic dynamics with detailed balance,
$\vb(\vx)=-\mD(\vx)\nabla\varphi^{ss}(\vx)$.  Therefore,
at the global minimum of $\varphi^{ss}(\vx,\vy;\Delta t)$, the flux $\dot{\vx}|_{\vx=\vx^*}=\vb(\vx^*)=0$.  This is the desired property
for an equilibrium state. In general, however, if without detailed balance we have $\vb(\vx)=-\mD(\vx)\nabla\varphi^{ss}(\vx)+\vg(\vx)$ where $\vg(\vx)\cdot\nabla\varphi^{ss}(\vx)=0$.
In this case, the global minimum of $\varphi^{ss}(\vx,\vy;\Delta t)$
implies $\nabla\varphi^{ss}(\vx^*)=0$ and a non vanishing flux
$\dot{\vx}|_{\vx=\vx^*}=\vg(\vx^*)$; it is a nonequilibrium steady state.

\begin{figure}[t]
\[
\includegraphics[scale=0.55]{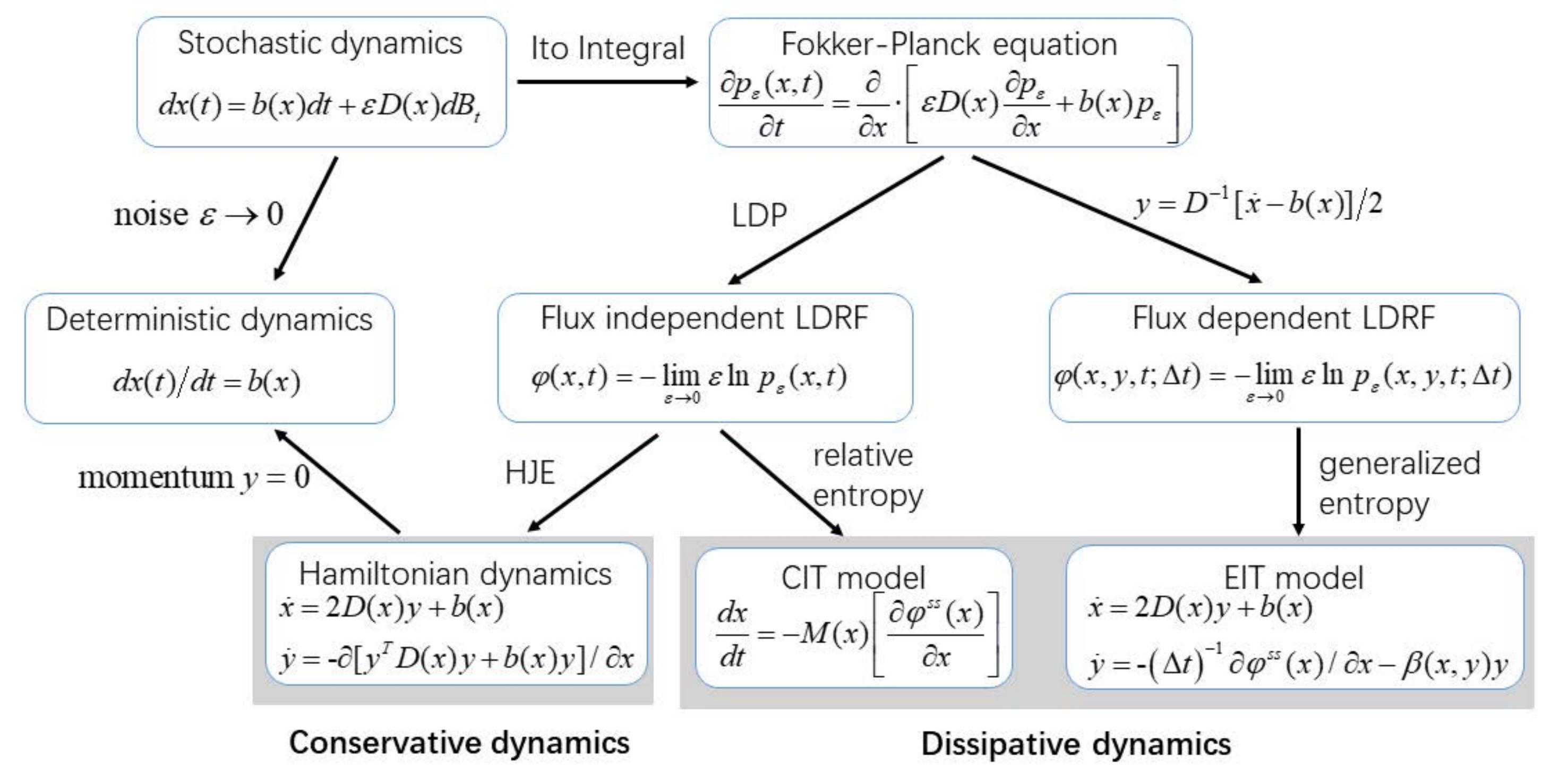}
\]
\caption{Diagram for stochastic dynamics, deterministic dynamics, Hamiltonian dynamics and dissipative dynamics by CIT and EIT, with large deviations rate function as a bridge in between. LDP, LDRF and HJE are short for large deviations principle, large deviations rate function and Hamiltonian-Jacobi equation respectively.}
\label{fig-2}
\end{figure}

\subsection{Explicit results for the Ornstein-Unlenbeck process}
Now we turn to an exactly solvable example -- the 1-d Ornstein-Uhlenbeck process (OUP):
$\rd x(t) = -bx(t)\rd t + \sqrt{2\epsilon D} \rd B(t)$,
where $b>0$, with its Kolmogorov forward equation (KFE) for the
transition probability $T_{\epsilon}(x,t|x',t')$,
\begin{equation}
   \frac{\partial T_{\epsilon}(x,t|x',t')}{\partial t} =
   \frac{\partial}{\partial x}\left( \epsilon D\frac{\partial T_{\epsilon}}{\partial x}
                + bx T_{\epsilon}\right), \  \
     T_{\epsilon}(x,t|x',t')|_{t=t'}=\delta(x-x').
\label{ou-fpe}
\end{equation}
Eq. \ref{ou-fpe} can be solved exactly to yield
\begin{equation}
      T_{\epsilon}(x,t|x',t') = \left\{\frac{b}{ 2\epsilon \pi D\big[1-e^{-2b(t-t')}\big]}
  \right\}^{\frac{1}{2}}\exp
  \left\{-\frac{ b\big[x-x' e^{-b(t-t')} \big]^2}{
         2\epsilon D\big[1-e^{-2b(t-t')}\big]} \right\}.
\end{equation}
More generally, Fokker-Planck equation (FPE) for the OUP
is the same linear partial differential equation in (\ref{ou-fpe})
with the $T_{\epsilon}$ replaced by a probability
density function
\begin{align}
p_{\epsilon}(x,t)=\bigg[\frac{b}{2\pi\epsilon D(1-e^{-2bt})}\bigg]^{1/2}\exp\bigg[-\frac{bx^2}{2\epsilon D(1-e^{-2bt})}\bigg],\  \
     p_{\epsilon}(x,0)=\delta(x),
\end{align}
that changes with time.

Based on these formulas, we can derive the flux-independent and flux-dependent large deviations rate function explicitly as
\begin{align}
&\varphi^{ss}(x)=\frac{bx^2}{2D},\\
&\varphi^{ss}(x,y;\Delta t)=\frac{bx^2}{2D}+(\Delta tD)y^2,
\end{align}
where $y=(\dot{x}+bx)/(2D)$. Repeating the same procedure of previous derivations, a natural dissipative dynamics suggested by EIT is
\begin{align}
&\frac{dx}{dt}=2Dy-bx,\\
&\frac{dy}{dt}=-\frac{bx}{\Delta t D}-\frac{by}{2},
\end{align}
by setting $\alpha(x,y)=0$. In this case, $\varphi^{ss}(x,y;\Delta t)$ turns to be the relative entropy with the dissipation rate as $(bx)^2/D+\Delta tbDy^2$.
Meanwhile, we can also get a Hamiltonian dynamics
\begin{align}
&\frac{dx}{dt}=2Dy-bx,\\
&\frac{dy}{dt}=by,
\end{align}
with the Hamiltonian function $H(x,y)=Dy^2-bxy$. It is noted that both dynamical systems are extensions of $dx/dt=-bx$, but their time reversibilities are completely opposite.

\subsection{Uncertainties in the zero-noise limit}

	What is the origin of the ``macroscopic, deterministic
thermodynamics?''  The title suggests an answer. This
seemingly paradoxical statement is precisely a consequence of
the concept of asymptotic limit, which had been considered as a
``devil's invention''.  Together with Zeno's paradox and Newton's
fluxion, they are a permanent part of the modern mathematics.  Furthermore, in theoretical physics, it is well appreciated that
when a limit process is singular, a wide range of counterintuitive
subjects can arise; and new theories of reality
emerge \cite{vulpiani-book}.

Let us again use the OUP to illustrate our idea. Consider
\begin{subequations}
\begin{eqnarray}
   p_{\sigma}(x,0) &=& \frac{1}{\sqrt{2\pi\sigma^2}}e^{-\frac{(x-x')^2}{2\sigma^2}},
\label{tozero}
\\
     p_{\sigma}^{\epsilon}(x,t) &=& \int_{\mathbb{R}} T_{\epsilon}(x,t|x'',0) p_{\sigma}(x'',0)\rd x''.
\end{eqnarray}
\end{subequations}
Noting the $p(x,0)$ in (\ref{tozero}) tending to $\delta(x-x')$
as $\sigma\to 0$.  We are particularly interested in the
limit of $\sigma\to 0$ and the ``zero-noise limit'' $\epsilon\to 0$.

	When considering WKB ansatz, we immediately notice that
the supposition $\delta(x-x')=e^{-\frac{1}{\epsilon}\varphi(x,0)}$
cannot be valid.
In other words, in the asymptotic limit, $\varphi(x,t)$
in terms of its characteristic lines is not fully defined by
$x=x'$ at $t=0$.  Additional information is required. This
additional information is precisely in the limit process of
$\sigma\to 0$.   On the other hand,
$p(x,0)=e^{-\varphi(x,0)/\epsilon}$
implies $\varphi(x,0) = -\epsilon\ln p_{\sigma}(x,0)$.
Therefore, in the limit of $\epsilon\to 0$, the $\varphi(x,0)$
corresponding to any proper $p_{\sigma}(x,0)$
vanishes.

	These uncertainty about $\varphi(x,0)$ is precisely solved
by the conjugate variable $y$ in the Hamiltonian
characteristic lines for the solution of the nonlinear HJE
\begin{equation}
  \frac{\partial\varphi(x,t)}{\partial t} = -D\left(
      \frac{\partial\varphi}{\partial x}\right)^2 + bx
    \left(\frac{\partial\varphi}{\partial x}\right).
\end{equation}
The ``momentum variable'' $y$ in the Hamiltonian
dynamics represents the randomness that gives
rise to a rare event in a stochastic dynamics.
Comparing the equation
\begin{equation}
  \frac{\rd x}{\rd t} = \frac{\partial H(x,y)}{\partial y}
           = -bx + 2Dy
\end{equation}
with the SDE
\begin{equation}
   \rd x(t) = - bxdt + \xi(t),\ \ \xi(t) =   \sqrt{2\epsilon D}\rd B(t),
\label{equation016}
\end{equation}
where $\xi(t)$ is a ``white noise'', we have
\begin{equation}
  y(t) = \sqrt{\frac{\epsilon}{2}}D^{-\frac{1}{2}}(x)
 \left(\frac{\rd B(t)}{\rd t}\right).
\label{equation017}
\end{equation}
Therefore, in terms of the white noise in (\ref{equation016}),
\begin{equation}
     y(t)\cdot \xi(t) =  \frac{ \epsilon }{\Delta t}.
\label{equation018}
\end{equation}
This is a kind of ``uncertainty principle'' between the
variance in $x$ and in momentum. Therefore, while $\varphi(x,t)$ emerges as a quantity
in the zero-noise limit, {\em its is neither the asymptotic
limit of the solution to FPE with proper initial value, nor
an asymptotic limit of the solution to KFE with
Dirac-$\delta$ initial value!}  The HJE
represents a novel behavior of its own.

	We now investigate the double limit
$\epsilon,\sigma\to 0$ for the function
\begin{subequations}
\label{equation14}
\begin{equation}
   -\epsilon\ln p^{\epsilon}_{\sigma}(x,t) =
       \frac{\epsilon\big(x-\mu(t)\big)^2}{2\Xi(t)}
       + \frac{\epsilon}{2}\ln\big( 2\pi\Xi(t) \big),
\end{equation}
in which, from Eq. \ref{tozero}, $\mu(t)=x'e^{-bt}$,
which is independent of $\epsilon$ and $\sigma^2$.
And,
\begin{equation}
	\theta^2(t) =  \frac{\epsilon D}{b}\Big(1-e^{-2bt}\Big),
\ \
    \Xi(t) = \sigma^2e^{-2bt}+\theta^2(t).
\end{equation}
\end{subequations}
The total Gaussian variance at time $t$, $\Xi(t)$,
has two parts, a decreasing contribution from the initial
$\sigma^2$ and an increasing Markovian $\theta^2(t)$.
In the limit of $\epsilon\to 0$ and $\sigma\to 0$,
\begin{eqnarray}
    && -\lim_{\sigma\to 0}\lim_{\epsilon\to 0} \epsilon\ln p^{\epsilon}_{\sigma}(x,t)
    = 0 \neq  -\lim_{\epsilon\to 0}\lim_{\sigma\to 0} \epsilon\ln p^{\epsilon}_{\sigma}(x,t)
\nonumber\\
\label{equation15}
	&=& \lim_{\epsilon\to 0}
         \left\{  \frac{b\big(x-x'e^{-bt}\big)^2}{2D\big(1-e^{-2bt}\big)}
       + \frac{\epsilon}{2}\ln\left[\frac{2\pi\epsilon D}{b}\Big(1-e^{-2bt}\Big) \right] \right\} = \frac{b\big(x-x'e^{-bt}\big)^2}{2D\big(1-e^{-2bt}\big)}.
\end{eqnarray}
The limit is highly singular; we particularly note that in the rhs
of (\ref{equation15}), there is an uncertainty at $t=0$, even
after taking the limit $\epsilon\to 0$.

%
%

\section{Discussion}

\subsection{Diffusion, friction, and mass}

{\bf\em The Einstein relation.} From a stochastic
treatment of mechanical motion, pioneered by Einstein,
Smoluchowski, and Langevin more than a century
ago, one has for example
\begin{equation}
   m\frac{\rd^2 x}{\rd t^2} = -\eta\frac{\rd x}{\rd t}
     - U'(x) + A\xi(t),
\end{equation}
respectively, in which $\xi(t)$ is a white noise represented
by the ``derivative'' of the non-differentiable Brownian motion,
$\rd B(t)/\rd t$.  Two limiting cases are particularly worth
discussion: (i) overdamped limit where $m=0$ and (ii) spatial
translational symmetric $U(x)=$ const.
The stationary distributions for (i) and (ii) are
\begin{equation}
     f_{x}(x) = Z_1^{-1} e^{-\frac{ 2\eta U(x) }{A^2}}
  \text{ and }
     f_{v}(v) = Z_2^{-1} e^{-\frac{m\eta \dot{x}^2}{A^2}}
\label{eqn35}
\end{equation}
in which $Z_1$ and $Z_2$ are corresponding
normalization factors for the two distributions.
Comparing (\ref{eqn35}) with Boltzmann's law and the
Maxwell distribution, one identifies $A^2=2\eta k_BT$,
where $k_B$ is Boltzmann's constant and $T$ is temperature
in Kelvin.  According to the diffusion theory,
$\frac{1}{2}(A/\eta)^2=D$ is the diffusion coefficient.
Therefore we arrive at the Einstein relation
$D\eta = k_BT$, a well known result in statistical mechanics.

{\bf\em Diffusion and mass.} In our present work, in the process of providing both
entropy in CIT, $-\varphi^{ss}(\vx)$, and flux-dependent
entropy in EIT, $-\varphi^{ss}(\vx,\vy;\Delta t)$, with a stochastic
dynamic foundation in a broad sense, we have been led to
an intriguing relation between the diffusion matrix $\mD(\vx)$
defined on the state space and the geometry concept of
an Riemannian metric in the tangent space for $\dot{\vx}$. The relation in (\ref{eqn023}) suggests an identification of
$k_BT[2\tau_{\delta}\mD(\vx)]^{-1}$ with a space-dependent
``mass'',  if $\vx$ is the Newtonian spatial coordinate.
Combining this with the Einstein relation,
$\frac{k_BT}{\eta}=\frac{(\Delta x)^2}{2(\Delta t)}=\frac{k_BT}{2\tau_{\delta}m}$.    This relation gives an provactive
hypothesis that $m\sim \frac{ (\Delta x)^{-2} }{k_BT}$.

\subsection{Fick's law as a consequence of Brownian motion}

	The heat or diffusion equation is obtained traditionally
by combining the continuity equation
$\partial u/\partial t = -\partial J/\partial x$ with
Fick's law $J=-D(\partial u/\partial x)$.  However, derivation
as such immediately suggests the possibility of generalizing
Fick's law.  But this turns out to be mis-leading.  In the
context of Brownian motion, the Fick's law should be
understood as ``an inbalance between the probability flux
$J_{A\to B}$ of a single diffusant, from region $A$ to region
$B$, and the $J_{B\to A}$.''  It is not driven
by concentration gradient {\it per se}; rather it is driven by
an ``entropic force'' $F$: $J = (F/\eta)u(x,t)$ where
$\eta$ is the frictional coefficient of the diffusant,
$F=-k_BT\partial\ln u(x,t)/\partial x$, and
$D=k_BT/\eta$ is the Einstein relation.
Any attempt to imporving Fick's law can only be considered as
a phenomological theory; a fundamental approach to the
subject has to consider {\em hydrodynamic limit} of
interacting particle systems \cite{gpv}.

\subsection{Parabolic vs. hyperbolic dynamics, and EIT}

	Another key anchoring points of EIT is the
{\em parabolic vs. hyperbolic dynamic equations}.
It is well-known that the former, in terms of diffusion, has
an infinite velocity for propagating a disturbance:  Solution
to $\partial u(x,t)/\partial t = \kappa\partial^2 u/\partial x^2$,
if $u(x,0)=\delta(x-x_0)$, $u(x,t)\neq 0$ for all $x\in\mathbb{R}$
when $t>0$.  This diffusive behavior is in sharp contrast to
hyperbolic dynamics.   Indeed, for many physical phenomena
on a short time scales and with high frequencies, inertia
plays an important role; the diffusive description becomes
unrealistic.  We would like to point out, however, that a more
fundamental distinction between parabolic vs. hyperbolic
dynamics is between {\em stochastic} and {\em deterministic}.
The latter emerges in a macroscopic limit.

\section*{Acknowledgements}
L.H. acknowledges the financial supports from the National Natural Science Foundation of China (Grants 21877070).


\begin{thebibliography}{00}

\bibitem{degroot}
de Groot, S. R. and Mazur, P. (1962)
{\em Non-Equilibrium Thermodynamics},
North-Holland, Amsterdam.

\bibitem{chapman-cowling}
Chapman, S. and Cowling, T. G. (1939)
{\em The Mathematical Theory of Non-Uniform
Gases},  Cambridge Univ. Press, U. K.

\bibitem{Muller}
M\"{u}ller, I., and Ruggeri, T. (1998)
{\em Rational Extended Thermodynamics}, Springer, New York.

\bibitem{jou}
Jou, D., Casas-V\'{a}zquez, J. and Lebon, G. (2009)
{\em Extended Irreversible Thermodynamics},
4$^{th}$ ed., Springer, New York.

\bibitem{gallavotti-book}
Gallavotti, G. (1999)
{\em  Statistical Mechanics: A Short Treatise},
Springer, Berlin.

\bibitem{qkkb}
Qian, H., Kjelstrup, S., Kolomeisky A. B. and Bedeaux D.
(2016) Entropy production in mesoscopic stochastic
thermodynamics: nonequilibrium kinetic cycles driven
by chemical potentials, temperatures, and mechanical
forces (topical review).
{\em J. Phys. Condens. Matter.} {\bf 28},
153004.

\bibitem{gq-pre-10}
Ge, H. and Qian, H. (2010)
The physical origins of entropy production, free energy
dissipation and their mathematical representations.
{\em Physical Review E}, {\bf 81}, 051133.

\bibitem{ge-qian-16}
Ge, H. and Qian, H. (2016)
Mesoscopic kinetic basis of macroscopic chemical thermodynamics: a
mathematical theory.
{\em Phys. Rev. E} {\bf 94}, 052150.

\bibitem{ge-qian-17}
Ge, H. and Qian, H. (2017)
Mathematical formalism of nonequilibrium thermodynamics for nonlinear
chemical reaction systems with general rate law.
{\em J. Stat. Phys.} {\bf 166}, 190--209.

\bibitem{hill-book}
Hill, T. L. (1977)
{\em Free Energy Transduction in Biology: The Steady-State
Kinetic and Thermodynamic Formalism},
Academic Press, New York.

\bibitem{wangjin}
Fang, X., Kruse, K., Lu, T. and Wang, J. (2019)
Nonequilibrium physics in biology.
{\em Rev. Mod. Phys.} to appear.

\bibitem{qian-arXiv}
Qian, H. (2017)
Kinematic basis of emergent energetic descriptions of
general stochastic dynamics.
arXiv:1704.01828.

\bibitem{qian-book-chapter}
Qian, H. (2019)
Nonlinear stochastic dynamics of complex systems, I.
In {\em Complexity Science: An Introduction},
Peletier, M. A., van Santen, R. A. and Steur, E. eds.,
World Scientific, Singapore, pp. 347--373.

\bibitem{qian-pla-2014}
Qian, H. (2014)
The zeroth law of thermodynamics and volume-preserving
conservative system in equilibrium with stochastic damping.
{\em Phys. Lett. A} {\bf 378}, 609--616.

\bibitem{ye-qian-19}
Ye, F. X.-F. and Qian, H. (2019)
Stochastic dynamics II: Finite random dynamical systems, linear representation, and entropy production.
{\em Discrete \& Continuous Dynamical Systems B}
{\bf 24}, 4341--4366.

\bibitem{vulpiani-book}
Chibbaro, S., Rondoni, L. and Vulpiani, A. (2014)
{\em Reductionism, Emergence and Levels of Reality},
Springer, New York.

\bibitem{gq-ijmpb}
Ge, H. and Qian, H. (2012)
Analytical mechanics in stochastic dynamics: Most probable
path, large-deviation rate function and Hamilton-Jacobi equation (review).
{\em International Journal of Modern Physics B},
{\bf 26}, 1230012.

\bibitem{zhu-15}
Zhu, Y., Hong, L., Yang, Z., and Yong, W. A.
(2015)
Conservation-dissipation formalism of irreversible thermodynamics.
{\em J. Non-Equil. Therm.} {\bf 40(2)}, 67--74.

\bibitem{gpv}
Guo, M. Z., Papanicolaou, G. C. and Varadhan, S. R. S. (1988)
Nonlinear diffusion limit for a system with nearest neighbor interactions.
{\em Commun. Math. Phys.} {\bf 118},  31--59.

\end{thebibliography}
\end{document}